\documentclass[10pt,a4paper]{article}
\usepackage{amssymb,amsmath,amsthm,amscd,enumerate}

\advance\textwidth 3cm
\advance\oddsidemargin -1,5cm
\advance\evensidemargin -2cm

\newcommand{\sym}{\mbox{Sym}_A(B\,\check{}) }
\newcommand{\gcp}{generic characteristic polynomial}
\title{Monogenous Algebras. Back to Kronecker.}
\author{Daniel Ferrand}
\begin{document}
\maketitle

\begin{center}
{\bf Introduction}
\end{center}

 In this note we develop some properties of those $A$-algebras
$B$ which can be generated by a single element after, if
need be, some faithfully flat base change $A \rightarrow A'$. Here they are
called \emph{locally simple}, instead of the less euphonious "monogenous". For
a finite algebra this condition is satisfied as soon as the geometric fibers
are single generator algebras. This property appears to be commonly
satisfied. In particular the morphisms between rings of algebraic
integers are locally simple.\par
A theme we insist
on expands an idea Kronecker introduced at the early beginning of the
algebraic theory of numbers, namely that many properties of a finite
free
$A$-algebra $B$ can be read through the characteristic polynomial of the
\emph{generic} element of $B$. To be precise, let $e_1, \ldots, e_n$ be a
basis of $B$ as an
$A$-module, and denote $S = A[T_1, \ldots, T_n]$. The generic element,
which may
be written as $\sum T_ie_i \in S\otimes_AB = B[T_1, \ldots, T_n]$, is a root of
its characteristic polynomial $F(X) \in S[X]$. Therefore we dispose of a
canonical morphism, called here the \emph{Kronecker morphism}
$$
S[X]/(F) \rightarrow S\otimes_AB.
$$
\hbox{We show that this morphism is
universally injective if and only if $B$ is locally simple over $A$.}\par
In the context of rings of algebraic integers the idea of Kronecker was
introduced and used by Hilbert
in his \emph{Zahlbericht}. We finally interpret some results of this
famous memoir from the point of view  previously introduced.

\vspace{5mm}

{\it In this note, all the rings are assumed to be commutative and to possess a
unit element, and all the ring morphisms are assumed to map unit element to
unit
element.}\\

\begin{center}
\textbf{1. Locally simple morphisms.}
\end{center}

\textbf{Definition 1.1}\,{\it A morphism $A \rightarrow B$ between rings is
called
\emph{simple} if
$B$ can be  generated, as an $A$-algebra, by a single element, in other
words if there
exists a surjective morphism of $A$-algebras $ A[X] \rightarrow B$.\par
A morphism $A
\rightarrow B$ is called
\emph{locally simple}  if there exists a faithfully flat morphism $ A
\rightarrow A'$ such that $A' \rightarrow A'\otimes_AB$ is simple.}\bigskip

Below we will give some specific properties of these morphisms. Let us
first recall how crucial they are to the theory of the norm functor (see
\cite{F}). With
 any finite and locally free morphism $A\rightarrow B$ of rank $d$ is
associated a covariant functor
$$
\textrm{N}_{B/A} : B-\textrm{{\bf Mod}} \; \longrightarrow\; A-\textrm{{\bf
Mod}},
$$
which extends the usual one defined for invertible $B$-modules $L$ (roughly
speaking,
by taking the norm of a cocycle of $L$). One has $\textrm{N}_{B/A}(B) = A$
but this functor is far from being additive. Let $F$ be a locally free
$B$-module of rank $n$. If $B$ is locally simple over
$A$ then  $\textrm{N}_{B/A}(F)$ is a locally free $A$-module of rank $n^d$.
But, in general, this $A$-module may have torsion,
even if
$B$ is a complete intersection over
$A$ (see
\cite{F} 4.3.4 and 4.4).
\vspace{1cm}

\textbf{Examples 1.2}\,  Consider a ring $A$ and the diagonal morphism $A
\rightarrow
A^n$. An element $x = (x_1,\dots, x_n) \in A^n$ is a generator of that
$A$-algebra if and only if the powers \, $1, x, x^2, \dots, x^{n-1}$ \,
form a basis of
the
$A$-module
$A^n$. Writing down these powers in the canonical basis of $A^n$, one sees
that $x$ is
a generator of the $A$-algebra $A^n$ if and only if the Van der Monde
determinant
$$ \prod_{i<j} (x_j-x_i)$$is invertible in $A$.\par
The existence of a sequence $(x_1,\dots, x_n)$ with this property is
clear if $A$
contains an infinite field. It is also clear that such a sequence cannot
exist if $A$
is too small; thus,
$\mathbb{F}_p
\rightarrow
\mathbb{F}_p^n$ is
\emph{not} simple if $n > p$. This remark, when $p = 2$, implies that
$\mathbb{Z}\rightarrow\mathbb{Z}^n$ is
\emph{not} simple if
$n
\geq 3$.
\par
On the other hand, there is a canonical way to adjoin to any ring $A$ a
sequence of
$n$ elements
$(x_1,\dots, x_n)$ making the Van der Monde determinant invertible. Just
take the ring
of fractions
$A' = A[X_1, \dots, X_n]_V$\, , where \, $V = \prod_{i<j}(X_j-X_i)$ and, for
$x_i$, take the
image in
$A'$ of $X_i$ ; the morphism $A \rightarrow A'$ is faithfully flat (and
smooth), and the
morphism $A' \rightarrow {A'}^n$ is simple ; thus for any $n$ and any ring
$A$, the morphism
$A \rightarrow A^n$ is locally simple.\par

A slight generalization implies that any finite \'etale morphism
$A\rightarrow B$ is locally
simple, because it is locally of the form $A \rightarrow A^n$.\medskip

\textbf{Example 1.3}\, More generally, let $A$ be a ring, and let $B_1,
\dots, B_s$ be a
sequence of finite and locally simple $A$-algebras. The product
$B_1\times \dots
\times B_s$ is locally simple over $A$.\par
 To see this, it is enough, by
induction on $s$,
to prove the result for two factors which we now denote by $B$ and $C$. Let
us choose generators
$b
\in B$ and  $c \in C$ and monic polynomials $P(T)$ and $Q(T)$ in $A[T]$
such that $P(b) =
0$ and $Q(c) = 0$ ; one then has a surjective morphism
$$
A[T]/(P)\times  A[T]/(Q)  \rightarrow B\times C ,
$$
and it is enough to prove that the product
$A[T]/(P)\times A[T]/(Q)$ is locally simple over $A$. Consider the
ring of
fractions $A' = A[X]_{R(X)}$ where we have made invertible the
\emph{resultant} (\cite{A} IV 6.6)
$$
R(X) = \mbox{res}_T(P(T+X), Q(T)).
$$
Let $x$ be the image of $X$ in $A'$. The standard property of the resultant
(see e.g \cite{A} IV 6.6 Cor.1 to Prop. 7), implies that the polynomials
$P(T+x)$ and $Q(T)$ are co-maximal in
$A'[T]$\, (i.e. they generate the unit
ideal). Therefore, the "Chinese remainder theorem" shows that the morphism
$$A'[T]
\longrightarrow A'[T]/(P(T+x)) \times A'[T]/(Q(T))$$
 is surjective. Moreover, the
$A'$-algebras $A'[T]/(P(T))$ and $A'[T]/(P(T+x))$ are clearly isomorphic.
Therefore, it
remains to show that the morphism $A \rightarrow A'$ is faithfully flat.
Since this morphism is clearly  flat
we have to show that any prime ideal $\mathfrak{p}$ of $A$ is the
restriction of a prime
ideal of $A'$. Let $A\rightarrow K$ be the morphism of $A$ to an algebraic
closure $K$
of the residue field $\kappa(\mathfrak{p})$; it is enough to see that this
morphism
factors through $A'$. Let us consider the images in $K[T]$ of the two monic
polynomials
$P(T)$ and $Q(T)$. Since $K$ is algebraically closed these polynomials
split, and we can translate the roots of
$P(T)$ away from those of
$Q(T)$. Thus, there exists
$x
\in K$ such that
$P(T+x)$ and
$Q(T)$ have no common root, i.e. such that the resultant $R(x)$ is non zero
in $K$. This element $x$ gives rise to the required morphism $A' =
A[X]_{R(X)}
\rightarrow K$.\bigskip

\textbf{Proposition 1.4} \textit {Let $B$ be a finite $A$-algebra.The following
conditions are equivalent:\\
i) The morphism $A \rightarrow B$ is locally simple.\\
ii) There exists a morphism $A\rightarrow A'$ such that $A'\rightarrow
A'\otimes_AB$ is simple,
and such that} Spec\textit{$(A')\rightarrow$}Spec\textit{$(A)$ surjective
(i.e the flatness of the base
change is superfluous).\\
iii) For any morphism $A \rightarrow K$ where $K$ is an
algebraically closed field, each local factor of $K\otimes_AB$ is a simple
$K$-algebra.\\
iv) For any prime ideal $\mathfrak{p}$ of $A$, there exists a finite extension
$\kappa(\mathfrak{p})\rightarrow k$ such that $k\otimes_AB$ is simple over
$k$.}\medskip

 Recall that a finite algebra $R$ over a field is the direct product of the
local rings
$R_{\mathfrak{m}}$, where $\mathfrak{m}$ runs through the (finite) set of
the maximal
ideals; these local rings  are called the
\emph{local factors} of
$R$.\\

 The ingredients used in the following proof all come from EGA IV, but, for
the convenience of the reader, I will give some details instead of
scattered references.\\

{\bf Lemma 1.4.1 } {\it Let $A \rightarrow B$ be a finite morphism. We suppose
an $A$-algebra $A\rightarrow E$ exists with the property that $E\otimes_AB$ is
simple over $E$.
Then, there exists a sub-$A$-algebra $F \subset E$} of finite type {\it
such that
$F\otimes_AB$ is simple over $F$.}\bigskip

Proof : Let  $x=\sum_{i=1}^{n} {x_i\otimes b_i} \in E\otimes_AB$ be a
generator as
$E$-algebra; the sub-$A$-algebra $E' = A[x_1, \dots, x_n] \subset E$ is of
finite type. Let us consider the morphism
$$
E'[X] \longrightarrow E'\otimes_AB,
$$
which maps $X$ to $x$; its cokernel $M$ is an $E'$-module of finite type, as
$E'\otimes_AB$ is, and we have $E\otimes_{E'}M = 0$. We shall enlarge $E'$
inside $E$ in order to get a finite type algebra $F$ such that
$F\otimes_{E'}M = 0$. By induction on the number of generators of $M$, (and
by looking at the
\emph{quotients} of
$M$) we are reduced to
the case where $M$ is monogenous, i.e where $M$ is isomorphic to a
quotient $E'/I$. The
hypothesis,
$E\otimes_{E'}M = 0$, reads then as $E = IE$, i.e as a relation: $1 =
\sum_{j=1}^m
a_jy_j$ with
$a_j
\in I$ and $y_j \in E$. This relation is already true in the $A$-algebra
of finite type
$E'[y_1, \dots, y_m]$.\bigskip

{\bf Lemma 1.4.2} {\it Let $\mathfrak{p}$ be a prime ideal in a ring $A$,
and let
$\kappa(\mathfrak{p}) \rightarrow k$ be a finite field extension. There exist
$t\in A- \mathfrak{p}$, a finite free morphism $A_t\rightarrow C$ and an
isomorphism
$\kappa(\mathfrak{p})\otimes_AC \; \tilde{\rightarrow} \; k$.}\medskip

Proof : We write $S = A - \mathfrak{p}$. By induction on the number of
generators of the
$\kappa(\mathfrak{p})$-algebra
$k$, we are reduced to proving the following.\par

Let $A_t \rightarrow C$ be a finite free morphism such that $k =
\kappa(\mathfrak{p}) \otimes_AC$ is a field, and let $k\rightarrow k'=k[x]$
be a
finite {\it simple} field extension. Then there exist $s \in S$ and a finite
free morphism $C_s \rightarrow C'$ such that $\kappa(\mathfrak{p})\otimes_AC'
\simeq k'$.\par
In fact, let
$F(X)\in S^{-1}C[X]$ be a monic polynomial whose image modulo
$\mathfrak{p}$ is the minimal polynomial of $x$ (such a polynomial $F$ exists
because the
morphism
$S^{-1}C \rightarrow S^{-1}C/\mathfrak{p}S^{-1}C \simeq k$ is surjective).
If $s \in S$
denotes the product of the denominators of the coefficients of $F$, one has
$F \in
C_s[X]$. The morphism
$$A_{st}\rightarrow C_s\rightarrow C'= C_s[X]/(F)$$
 is then free, and one gets an
isomorphism
$\kappa(\mathfrak{p})\otimes_AC'\simeq k'$.\bigskip

Proof of the proposition. It is clear
that {\it i)} implies {\it ii)}. Let us prove that {\it ii)} implies {\it
iii)}.  Let
$A'$ be an
$A$-algebra such that
$A'\otimes_AB$ is generated by one element
and such that the map Spec$(A') \rightarrow \mbox{Spec}(A)$ is surjective.
By the above
lemma 1.4.1 there exists a sub-$A$-algebra $F\subset A'$, of finite type,
such that
$F\otimes_AB$ is simple over $F$. Let $A \rightarrow K$ be a morphism where
$K$ is an
algebraically closed field, and denote by
$\mathfrak{p}$ its kernel. By hypothesis, the prime ideal
$\mathfrak{p}$ is the restriction to
$A$ of a prime ideal $\mathfrak{p'}$ of $A'$; it is also the restriction
of the prime ideal $\mathfrak{q} = \mathfrak{p'} \cap F$ of $F$, therefore
$\kappa(\mathfrak{p})\otimes_AF \neq 0$. Then, as
$K$ is algebraically closed, the "Hilbert Nullstellensatz" (\cite{AC} V 3.3
Prop.1)
implies that the given morphism
$\kappa(\mathfrak{p})\rightarrow K$ factors through
$\kappa(\mathfrak{p})\otimes_AF$, i.e. that $A\rightarrow K$ factors through
$F$.
$$
\begin{CD}
A @>>> F @>>> F\otimes_AB\\
@VVV @VVV @VVV\\
\kappa(\mathfrak{p}) @>>> K @>>> K\otimes_AB
\end{CD}
$$
But the morphism
$F
\longrightarrow F\otimes_AB$ is simple.
Therefore, the $K$-algebra
$K\otimes_AB$ is simple, and a fortiori each of its factors is.\par

{\it iii)} $\Rightarrow$ {\it iv)}. Let $K$ be an algebraic closure of a
residue field
$\kappa(\mathfrak{p})$ of $A$. By the hypothesis {\it iii)} and the example
{\bf 1.3},
the $K$-algebra $K\otimes _AB$ is simple; by lemma 1.4.1, there exists
a finite
sub-extension $k \subset K$ such that $k\otimes_A B$ is a simple
$k$-algebra.\par

{\it iv)} $\Rightarrow$ {\it i).} Suppose first we have already shown that for
each prime ideal
$\mathfrak{p}$ of
$A$ there exist an element $t\in A-\mathfrak{p}$ and a finite free
morphism
$A_t\rightarrow C$ such that
$C\rightarrow C\otimes_AB$ is simple.\par
 Then the image of the morphism
Spec$(C)\rightarrow$
Spec$(A)$ is the open set $D(t)$, and it contains $\mathfrak{p}$. As
Spec$(A)$ is
quasi-compact, a finite number of such morphisms $A
\rightarrow C_i, i=1,\dots n$, are enough for covering Spec$(A)$. Hence we can
take $A' = C_1\times \dots \times C_n$; it is faithfully flat over $A$, and
$A'\rightarrow A'\otimes_AB$ is simple.\par
It remains to prove the existence of those required morphisms $A_t\rightarrow
C$. So let
$\mathfrak{p}$ be a prime ideal in
$A$. According to {\it iv)}, there exists a
finite extension $\kappa(\mathfrak{p}) \rightarrow k$
 such that
$k\rightarrow k\otimes_AB$ is simple. By lemma 1.4.2, one can choose a
$t \in
S = A -\mathfrak{p}$, a finite free morphism $A_t\rightarrow C$ and an
isomorphism $\kappa(\mathfrak{p})\otimes_AC \;\tilde{\rightarrow}\; k$. The
morphism
$C \rightarrow
\kappa(\mathfrak{p})\otimes_AC \simeq k$ is the composite of the surjection
$S^{-1}C
\rightarrow S^{-1}(C/\mathfrak{p}C)$ and of the localization $C \rightarrow \
S^{-1}C$. Then, a generator $\xi$ of $k\otimes_AB =
S^{-1}(C/\mathfrak{p}C)\otimes_AB$ may be lifted as an element $x \in
S^{-1}(C\otimes_AB)$. It is a generator of the $S^{-1}C\,$-algebra
$S^{-1}(C\otimes_AB)$.
$$
\begin{CD}
S^{-1}C @>>> S^{-1}C[x] @>>> S^{-1}C\otimes_AB\\
@VVV @VVV @VVV\\
k @>>> k[\xi] @= k\otimes_AB
\end{CD}
$$
In fact the cokernel
of the injective map
$S^{-1}C[x] \hookrightarrow S^{-1}(C\otimes_AB)$, is a finitely generated
module
over $ S^{-1}A = A _{\mathfrak{p}}\,$, which is zero modulo $\mathfrak{p}$. The
Nakayama lemma implies this cokernel is zero.\par
 Finally, there is a
$s'
\in S$ such that
 $x \in C_{s'}\otimes_AB$. Using again the above finiteness property of the
cokernel, we can find a $s'' \in S$  such that the map
$C_{s's''}[x]\;
\rightarrow
\; C_{s's''}\otimes_AB$ is an isomorphism. The morphism $A_{s's''t} \rightarrow
C_{s's''}$ has the
required properties.
\bigskip

\textbf{Corollary 1.5} \textit{A finite $A$-algebra $B$ is locally simple
if and only if
$\; \Omega_{B/A}^2 = 0.$}\medskip

Proof. If $B$ is simple over $A$, then the $B$-module $\Omega_{B/A}^1$ is
generated by
one element, namely the differential $\, d_{B/A}(x)\, $ of a generator $x$.
Therefore its square wedge is zero.
The same conclusion is true if $B$ is locally simple because of the isomorphism
$A'\otimes_A\Omega_{B/A}^1 \simeq \Omega_{A'\otimes_AB/A'}^1$.\par
 Conversely, suppose that
$\Omega_{B/A}^2 = 0$. We shall prove that the condition $iii)$ of {\bf
1.4} is satisfied. So let
$A\rightarrow K$ be a morphism to an algebraically closed
field $K$. Let $R$ be a local factor of $K\otimes_AB$. By assumption, one has
$\Omega_{R/K}^2 = 0$. We write  $\Omega = \Omega_{R/K}^1$, and we denote by
$\mathfrak{m}$ the maximal ideal of $R$. Since
$\wedge^2(\Omega/\mathfrak{m}\Omega) = 0$ the dimension of the
$R/\mathfrak{m}$-vector space $\Omega/\mathfrak{m}\Omega$ is $\leq 1$. As
$K$ is
algebraically closed,
$K \rightarrow R/\mathfrak{m}$ is an isomorphism. We will
use now the  well-known (see below)
$K$-linear isomorphism
$$
\delta :\mathfrak{m}/\mathfrak{m}^2 \quad \widetilde{\longrightarrow} \quad
\Omega/\mathfrak{m}\Omega .
$$
It implies that $\mathfrak{m}/\mathfrak{m}^2$ is a
$K$-vector space of dimension $\leq 1$. From the Nakayama lemma we then deduce
that the ideal
$\mathfrak{m}$ may be generated by one element. Thus $R$ is a simple
$K$-algebra.\par
(For lack of an elementary reference, we briefly recall that
\ $\delta$ \ is induced by the  differential \hbox{$d_{R/K} :
\mathfrak{m} \,
\rightarrow \, \Omega$,} and that the inverse of $\delta$ is defined as
follows. Let
$s: R\rightarrow R/\mathfrak{m} \simeq K$ be the canonical morphism. The
map
$ R
\rightarrow
\mathfrak{m}/\mathfrak{m}^2, \; x\mapsto$ class of
$x-s(x)$ mod.$\mathfrak{m}^2$, \; is a derivation. By the universal property
of
$\Omega$, this derivation extends to a linear map
$\Omega/\mathfrak{m}\Omega \rightarrow
\mathfrak{m}/\mathfrak{m}^2$, which is easily seen to be the inverse of
$\delta$.)
\bigskip

\textbf{Corollary 1.6} {\it Let $A \stackrel{u}{\rightarrow} B
\stackrel{v}{\rightarrow}
C$ be finite morphisms. Then the composite $vu$ is locally simple if either: \\
- $u$ is locally simple and $v$ is net (i.e unramified).\\
- $u$ is net and $v$ is locally simple.}\bigskip

This result, which generalizes {\bf 1.3}, is easily deduced from the previous
corollary and from
the exact sequence
$$
\Omega_{B/A}^1\otimes_AC  \rightarrow \Omega_{C/A}^1 \rightarrow \Omega_{C/B}^1
\rightarrow 0 .
$$
\bigskip

\textbf{Corollary 1.7} \textit{ Let $A$ be a Dedekind domain, $K\rightarrow
L$ a finite
separable extension of its field of fractions, and let $B$ be the integral
closure of
$A$ in $L$. Suppose that all the residue field extensions are separable. Then
$A\rightarrow B$ is locally simple.}\bigskip

Proof. Let $\mathfrak{n}$ be a maximal ideal of $B$, and let $\mathfrak{m}
= A \cap
\mathfrak{n}$. As $B_{\mathfrak{n}}$ is a discrete valuation ring, the $\kappa
(\mathfrak{n})$-vector space $\mathfrak{n}/\mathfrak{n}^2$ is of dimension
1. Since
$\kappa(\mathfrak{n})$ is supposed to be separable over
$\kappa(\mathfrak{m})$ one has \,
$\Omega_{\kappa(\mathfrak{n})/\kappa(\mathfrak{m})}^1 = 0$. Therefore, the
exact sequence
$$
\mathfrak{n}/\mathfrak{n}^2 \rightarrow \Omega_{B/A}^1\otimes_B
B/\mathfrak{n} \rightarrow
\Omega_{\kappa(\mathfrak{n})/\kappa(\mathfrak{m})}^1 \rightarrow 0
$$
shows that $\Omega_{B/A}^1\otimes_B B/\mathfrak{n}$ is a vector space of
rank $\leq 1$. Hence for each maximal ideal $\mathfrak{n}$ one has \;
$\Omega_{B/A}^2\otimes_B
B/\mathfrak{n} = 0$, and the Nakayama lemma gives
$(\Omega_{B/A}^2)_\mathfrak{n}
= 0$. Since this is true for each maximal ideal of $B$, we may conclude that
$\Omega_{B/A}^2 = 0$.
\newpage

\begin{center}
\textbf{2. The Kronecker morphism.}
\end{center}\bigskip

\textbf{2.1. The generic element.} Let $A \rightarrow B$ be a finite and
locally
free morphism. Denoting by $B\,\check{} = $Hom$_A(B, A)$
the dual of $B$, we let
$$ \xi \in B\,\check{}\otimes_A B $$
be the element corresponding to the identity map of $B$ via the isomorphism
$$ B\,\check{}\otimes_A B \quad \widetilde{\longrightarrow} \quad
\mbox{End}_A(B)$$
which sends $\beta\otimes x \in B\,\check{}\otimes_A B $  to
the map $b \mapsto
\beta(b)x$. If we write $\xi = \sum \beta_i\otimes b_i$, then, for all $b
\in B$, one has
$ b =
\sum \beta_i(b)b_i.$\par
 When viewing it as an element of
$\mbox{Sym}_A(B\,\check{}) \otimes_A
B$, we
call $\xi$ the \emph{generic element} of $B$, and we call
$\mbox{Sym}_A(B\,\check{})$
 the \emph{ring of parameters} for the elements of $B$. In
fact, an
element $x$ in $B$ uniquely determines the $A$-linear map $ B\,\check{}
\rightarrow
A$ given by $ u \mapsto u(x)$. This map extends to a morphism of $A$-algebras
$$
\gamma_x : \mbox{Sym}_A(B\,\check{}) \rightarrow A .
$$
The morphism $\gamma_x$ has to be seen as the \emph{specialization of
parameters}
associated with
$x$ because  we recover
$x$ as the image of the generic element
$\xi$ by the morphism
$$\gamma_x \otimes 1: \mbox{Sym}_A(B\,\check{})\otimes_A B \rightarrow B
.$$\par

If $(e_i)$ is a basis of $B$ (as $A$-module), and if $(e_i\check{})$
denotes the
dual
basis, one has : $\xi = \sum_{i} e_i\check{}\,\otimes e_i$. The ring of
parameters $\mbox{Sym}_A(B\ \check{})$ is then isomorphic to the polynomial
ring
$A[T_1, ... , T_n]$, where $T_i$ stands for $e_i\check{}$ , and the generic
element
is usually written
$$\xi = \sum_{i} T_i e_i .$$

But introducing variables may hide the important fact that the ring of
parameters
for $B$ is \emph{contra}variant in $B$. Understanding the functoriality of
$\xi$ is
thus easier if one keeps its intrinsic definition and uses the following
remark.\par

Let $u : B \rightarrow C$ be a morphism of $A$-algebras. Suppose $C$ to be
finite and
locally free over $A$. Under the canonical maps
$$
\mbox{Hom}_A(B, B) \longrightarrow \mbox{Hom}_A(B, C) \longleftarrow
\mbox{Hom}_A(C, C)
$$
the images of id$_B$ and of id$_C$ are both equal to $u \in \mbox{Hom}_A(B,
C)$. In
order to extend this to the generic elements let us consider the morphisms of
$A$-algebras
$$
 \mbox{Sym}_A(B\,\check{})\otimes_AB \stackrel{1\otimes u}{\longrightarrow}
 \mbox{Sym}_A(B\,\check{})\otimes_AC \stackrel{v\otimes 1}{\longleftarrow}
\mbox{Sym}_A(C\,\check{})\otimes_AC,
$$
where $v = \mbox{Sym}_A(u\,\check{})$. Then, the images of the generic elements
$\xi_B
\in \mbox{Sym}_A(B\,\check{})\otimes_AB$ and $\xi_C
\in \mbox{Sym}_A(C\,\check{})\otimes_AC$ are equal in the ring
$\mbox{Sym}_A(B\,\check{})\otimes_AC$.

\bigskip

\textbf{2.2 The Kronecker morphism.}
Let again $A \rightarrow B$ be a finite and locally free morphism. Let
$$
F_{B/A}(X)\in \sym[X]
$$
be the characteristic polynomial of the generic
element of
$B$. From now on this polynomial will be called the \emph{\gcp}.\par
The relation $
F_{B/A}(X) = 0$ is called by Hilbert ({\it Zahlbericht}, ch.IV,
\S 10) the
\emph{fundamental equation} of the $A$-algebra $B$. The
generic element is a root of this equation (Hamilton-Cayley theorem).
 Therefore there
exists a morphism of $\sym$-algebras
$$\sym[X]/(F_{B/A}) \longrightarrow  \sym\otimes_A B,$$
which maps (the class of) $X$ to $\xi$; it will be called the
\emph{Kronecker morphism} of $B/A$.\\

{\bf 2.2.1}\;Suppose  $B = A^n$, and choose the canonical basis $(e_i)$ for
$A^n$. The ring of parameters $\sym$ is then isomorphic to $S = A[T_1, \dots,
T_n]$, where
$T_i$ stands for
the
$i$-th projection $A^n \rightarrow A$. An immediate calculation gives
$$
F_{B/A}(X) = \prod_{i=1}^n (X-T_i),
$$
and the Kronecker morphism
$$
S[X]/(\prod (X-T_i)) \longrightarrow S^n
$$
is defined by $ X \mapsto (T_1, \dots, T_n)$. It is injective since the
Van der
Monde determinant is a regular element in $S$.\par
The coefficients of the \gcp \;  $F_{B/A}$\, are symmetric polynomials in
the $T_i$, i.e
they are invariant under the automorphisms of the $A$-algebra $B = A^n$.
Below we
will show this to be a general fact.\bigskip

{\bf 2.2.2}\;Direct calculations are seldom illuminating, even in the simplest
cases. Let, for example,
$B = A[Y]/(G)$ be the
$A$-algebra of rank 3 defined by the polynomial
$$
G(Y) = Y^3+a_2Y^2+a_1Y+a_0 .
$$
If we write the generic element of $B$ as $\xi = T_0 + T_1y + T_2y^2$, then
$$
F_{B/A}(X) = (X-T_0)^3 + (X-T_0)^2[a_2T_1+(2a_1-a_2^2)T_2] \hspace{6cm}$$
$$\hspace{2cm}+(X-T_0)[a_1T_1^2 + (3a_0-a_1a_2)T_1T_2 +(a_1^2-2a_0a_2)T_2^2]$$
$$
 \hspace{5cm}+[a_0T_1^3 - a_0a_2T_1^2T_2 +
a_0a_1T_1T_2^2 -a_0^2T_2^3] .
$$
From this formula, it is not even clear if the Kronecker morphism is
injective. In fact
it is (cf {\bf 2.4}).
\bigskip

{\bf 2.2.3}\; Let $B = A[u, v]$ with $u^2 = v^2 = 0$. It is a radicial
$A$-algebra of
rank 4. Writing the generic element as $\xi = T_0+T_1u+T_2v+T_3uv$, we find
$$
F_{B/A}(X) = (X-T_0)^4.
$$
Since $(\xi-T_0)^3 = 0$, the Kronecker morphism is \emph{not} injective in
that case.

\vspace{1 cm}

{\bf 2.3}\; A few words now on the functoriality of these notions.\par
 Let
$u:B
\rightarrow C$ be a morphism between finite and locally free $A$-algebras.
Denote the
rings of parameters by $S =
\sym$, and $T = \mbox{Sym}_A(C\,\check{})$, and let
$$
v = \mbox{Sym}_A(u\,\check{}) : T
\longrightarrow S
$$
be the morphism associated to $u$. Finally, let the norm maps relative
to $C/A$ be
denoted by $\textrm{N}_T = \textrm{N}_{T\otimes_AC/T}$, and $\textrm{N}_S =
\textrm{N}_{S\otimes_AC/S}$. We have the following commutative diagram
$$
\begin{CD}
S\otimes_AB @>{1\otimes u}>>S\otimes_AC@<{v\otimes 1}<< T\otimes_AC\\
&&@V{\small{\textrm{N}_S }}VV @VV{\small{\textrm{N}_T}}V\\
&&S @<<v<T
\end{CD}
$$
Let $F_{B/A}(X)\in S[X]$, and
$F_{C/A}(X) \in
T[X]$ be the generic characteristic polynomials. Since $1\otimes u(\xi_B) =
v\otimes
1(\xi_C)$, we see that
$$v(F_{C/A}) = \textrm{N}_S(X- v\otimes 1(\xi_C)) = \textrm{N}_S(X- 1\otimes
u(\xi_B)).
$$
We can make this last polynomial explicit in some particular cases.\\

{\bf 2.3.1}\; $C = B$ and $u$ is an automorphism of the $A$-algebra $B$.\par
As the norm commutes with any automorphism, we have $v(F_{B/A}) = F_{B/A}$.
Hence\\

{\it The \gcp\; } $F_{B/A} \in \sym [X]$\; {\it is invariant under any
automorphism of
the
$A$-algebra
$B$ acting (contra-variantly) on \, }$\sym$.\bigskip

{\bf 2.3.2}\; Suppose that the morphism $u:B\rightarrow C$ is locally free of
rank $d$. Then
by transitivity of the norm, we have
$$
\textrm{N}_{C/A} = \textrm{N}_{B/A}\circ \textrm{N}_{C/B}.
$$
Moreover, for $b\in B, \; \textrm{N}_{C/B}(u(b)) = b^d$. Therefore,
$$
 \textrm{N}_S(X- 1\otimes
u(\xi_B))
=\textrm{N}_{S\otimes_AB/S}(\,\textrm{N}_{S\otimes_AC/S\otimes_AB}(X-
1\otimes
u(\xi_B)) = F_{B/A}(X)^d,
$$
and we get
$$
v(F_{C/A}(X)) =  F_{B/A}(X)^d.
$$
In that case, the morphism $v : \textrm{Sym}_A(C\,\check{}) \rightarrow
\sym$ is
surjective and it may be seen as "sending to 0" the variables relative to
the quotient
$C/B$.\bigskip

{\bf 2.3.3}\; Consider the case where $C = B/J$ is a quotient by a nilpotent
ideal $J$ such that
the successive quotients $J^s/J^{s+1}$ are locally free $C$-modules of
constant rank;
let $e$ be the sum of these ranks. According to (\cite{A} III 9.4  Prop. 5),
we have in
$S[X]$
$$
F_{B/A}(X) = \textrm{N}_S(X- 1\otimes u(\xi_B))^{e},
$$
whence the equality
$$
F_{B/A}(X) = v(F_{C/A}(X)^{e}).
$$\bigskip

{\bf 2.3.4}\; Let us give a precise meaning to the intuitively
clear following sentence:\\
 {\it The
\gcp\; of a product of rings is the product of the generic characteristic
polynomials of the
factors.}\par
Let $B = \prod B_i$ be a decomposition in a finite product of $A$-algebras.
Denote by
\hbox{$p_i : B \rightarrow B_i$} the projections, and let
$$q_i=\textrm{Sym}(p_i\check{}):  S_i =\textrm{Sym}_A(B_i\check{})
\quad \longrightarrow \quad S = \textrm{Sym}_A(B\,\check{})$$
be the injective morphism associated to $p_i$. The decomposition
$S\otimes_AB \simeq
\prod S\otimes_AB_i$ shows that
$$
F_{B/A}(X) = \prod \textrm{N}_{S\otimes_AB_i/S}(X- 1_S\otimes p_i(\xi_B)),
$$
If we now apply the remark above with
$p_i : B\rightarrow B_i$ instead of $u:B\rightarrow C$, we get
$$
\textrm{N}_{S\otimes_AB_i/S}(X- 1_S\otimes p_i(\xi_B)) = q_i(F_{B_i/A}(X)).
$$
Putting these equalities together, we find the expected formula for the
product :
$$
F_{B/A}(X) = \prod  q_i(F_{B_i/A}(X)).
$$
\vspace{0,5cm}

\textbf{Theorem 2.4}\; (Injectivity of the Kronecker morphism) \textit{Let
$A \rightarrow
B$ be a finite and locally free morphism. Then the following conditions are
equivalent:}\par 
{\it i)  $B$ is locally simple over
$A$.}\par
{\it ii) The Kronecker morphism}
$$\sym[X]/(F_{B/A}) \longrightarrow \sym\otimes_A B,$$\textit{is injective,
and remains
injective after any base change $A\rightarrow A'$.}\medskip

Let us show the implication $i) \Rightarrow ii)$. We can clearly suppose
$B$ to be
simple, hence of the form
$A[Y]/(G)$, where $G$ is a monic polynomial of degree $n$. We write
$y$  for the class of  $Y$ in $B$, and we choose the basis $(1, y,\dots
,y^{n-1})$ for
$B$. The ring of parameters $\sym$ will then be seen as the polynomial ring
$S = A[T_0,
T_1,\dots ,T_{n-1}]$, in such a way that the generic element would be
written as
$$\xi = T_0+T_1y+\dots +T_{n-1}y^{n-1}.$$
Checking the injectivity of the Kronecker morphism amounts to proving the
following: any relation of the form
$$s_0 + s_1\xi + ... + s_{n-1}\xi^{n-1} = 0$$
with the $s_i$ in $S$, implies that all the $s_i$ are zero; in other words, one
has to show that the family  $(1, \xi, ...,\xi^{n-1})$ of elements of
$S\otimes_AB$ is free over $S$. For doing so, we consider the determinant
of the
matrix of the $\xi^{j}$ on the basis $(y^{i})$, and we show it is a regular
(i.e not
a zero divisor) element of $S$.\\
Let $U_{ij}\in S$ be the polynomials defined by
$$\xi^{j} = U_{0,j}+U_{1,j}y+\dots +U_{n-1,j}y^{n-1},$$
and let $U =\mbox{det}(U_{ij}).$\medskip

An explicit example will perhaps give the flavor of the situation.\par
\noindent When
$G(Y) = Y^3+a_2Y^2+a_1Y+a_0$ some by hand calculations give
$$U = T_1^{3} -
2a_2T_1^2T_{2}^{}+(a_1+a_2^{2})T_1^{}T_2^2+(a_0-a_1a_2)T_2^{3}.$$
In this example, a fact is to be noticed: $U$
is a monic polynomial in $T_1$. We will check this to be a general
fact.\medskip

Each of the polynomials $U_{ij}$ is homogeneous in $T_0,
T_1,\dots ,T_{n-1}$, of degree $j$. Therefore the determinant $U$ is a
homogeneous polynomial of degree $N = 1+2+ \cdots +n-1$. On the other hand,
$U(0, T_1, 0,\dots ,
0) = T_1^N$. In fact if we map $\xi$ to $T_1y$ then $\xi^j$ is mapped to
$T_1^jy^j$,
and the matrix involved becomes the diagonal matrix \; diag$(1, T_1^{}, \dots,
T_1^{n-1})$. These two facts together imply that
$U$ is a monic polynomial in $T_1$. Hence $U$ is a regular element in
$S$, and it remains regular after any base change $A\rightarrow A'$.\medskip

To prove the implication $ii) \Rightarrow i)$ let us first simplify the
notations and
write the Kronecker morphism as
$$
u : S[X]/(F) \longrightarrow S\otimes_AB.
$$
The cokernel $M$ of $u$ is an $S$-module of finite presentation, therefore
the set $V
\subset \textrm{Spec}(S)$ of those prime ideals $\mathfrak{n}$ of $S$ such that
$M_{\mathfrak{n}} = 0$ is open and quasi-compact (Hence $V$ may be
covered by a finite family $(V_i)$ of affine open subsets, and someone who
dislikes schemes could replace $V$ by the ring $A' = \prod_i \Gamma (V_i)$).
The injectivity assumption implies that
$\mathfrak{n}$ is in
$V$ if and only if
$u_{\mathfrak{n}}$ is an isomorphism. Thus for such an \ $\mathfrak{n}$ \  the
$S_{\mathfrak{n}}$-algebra
$S_{\mathfrak{n}}\otimes_AB$ is simple. As the
morphism
$V
\rightarrow
\mbox{Spec}(A)$ is clearly flat, we have only to check that it is
surjective. So, let
$k = \kappa(\mathfrak{p})$ be the residue field at a prime ideal
$\mathfrak{p}$ of $A$. The ring $k\otimes_AS$ is isomorphic to a polynomial
ring
over $k$; thus it is an integral
domain, and we denote its fraction field by $K$. By assumption,
the
$k\otimes_AS$-linear map $1_k\otimes u$ is injective; by localization, the
$K$-linear
map $u_K : K[X]/(F) \rightarrow K\otimes_AB$ is still injective. But both
sides are
$K$-vector spaces of the same dimension; therefore $u_K$ is bijective. Hence
the generic
point of $k\otimes_AS$ is a point of $V$, and it lies over $\mathfrak{p}$.

\vspace{1cm}

\textbf{Remark 2.5}\; An other proof of the implication $ii) \Rightarrow
i)$ uses
the condition $iii)$ of the proposition {\bf 1.4}. We will now give its
main step
because it
seems to be of interest in itself.\par
{\it Let $K$ be an algebraically closed field, and $R$ a
finite local $K$-algebra. We suppose that there exist a non zero
$K$-algebra $S$, a
monic polynomial $F(X) \in S[X]$ of degree} $n = \mbox{rank}_K(R)$,{\it and an
\emph{injective} morphism of $S$-algebras $u : S[X]/(F) \rightarrow
S\otimes_KR$. Then
$R$ is a simple $K$-algebra.}\par
 Proof : We write $R = K + J$ where $J$ is nilpotent. Let $m$ be the lowest
integer such
that $J^m = 0$; hence, in the filtration
$$ R \supset J \supset J^2 \supset ... \supset J^{m-1} \supset J^m = 0$$
all those  $K$-subspaces are distinct. Therefore, we have $m \leq \dim_K(R) =
n$. Let
$x$ denote
the class of $X$ in $S[X]/(F)$. We write  $u(x) = s +
\eta \in S\otimes_KR  = S + S\otimes_K J$, with
$s
\in S$ and
$ \eta \in S\otimes_A J$. Since $u((x-s)^m) = \eta^m = 0$, the injectivity of
$u$ implies
that $F(X)$ divides  $(X-s)^m$. Therefore $m = n$ because $\deg(F) = n \geq
m$. Thus,
$J^{n-1} = J^{m-1} \neq 0$. But $J$ is a vector space of dimension $n-1$, and
the filtration above is strict; therefore the vector space $J/J^2$ is of
rank one, i.e the
ideal $J$ is generated by one element (Nakayama) and we conclude that $R$ is
simple.

\vspace{0,5 cm}

The following consequence has most probably been already noticed, at least
in its polynomial setting.\bigskip

\textbf{Corollary  2.6} \textit{Let $A \rightarrow B$ be a finite and
locally free
morphism. Suppose $A$ and $B$ to be domains. If $B$ is locally simple over
$A$, then the
norm of the generic element of $B$ generates a prime ideal in } $\sym$. {\it In
particular, let
$L/K$ be a finite simple field extension; choose any $K$-basis $(e_1,\dots
,e_n)$ for
$L$. Then the polynomial}
$$F(T_1,\dots ,T_n)=\mbox{Norm}_{L/K}(T_1e_1+\dots +T_ne_n)$$
\textit{is irreducible in $K[T_1,\dots ,T_n]$.}\\

Proof : We write $\sym = S$, and we consider the morphism of $A$-algebras
\hbox{$ u : S
\rightarrow S[X]$} extending the linear map $B\,\check{} \rightarrow S[X]$
given by $
\beta \mapsto \beta(1)X - \beta$. If we denote by $v :S[X] \rightarrow S$ the
morphism of $S$-algebras which sends $X$ to $0$, then $v\circ u$ is clearly an
automorphism of $S$. Let $\xi \in S\otimes_AB$ be the generic element of
$B$. Its
image $u\otimes1(\xi)$ in $S[X]\otimes_AB$ is easily seen to be $X- \xi$. The
commutativity of the
square
$$
\begin{CD}
S\otimes_AB @>{u\otimes 1}>>S[X]\otimes_AB\\
@V{\small{\mbox{Norm}}}VV @VV{\small{\mbox{Norm}}}V\\
S @>>u>S[X]
\end{CD}
$$
shows that $u(\mbox{Norm}_{B/A}(\xi)) = F_{B/A}(X)$. Therefore, $u$ induces a
morphism
$$
\bar{u} : S/\mbox{Norm}(\xi)S \longrightarrow S[X]/(F).
$$
Using the morphism $v$, we see that $\bar{u}$ is injective. The conclusion
now follows
from the implication $i) \Rightarrow ii)$ of the theorem, and from the fact
that
$S\otimes_AB$ is a domain.
\bigskip

\textbf{Remark  2.7} The simplest non simple (!!) field extension is
$$ K = \mathbb{F}_2(X, Y) \; \subset \; L =  \mathbb{F}_2(U, V),$$
given by $X = U^2,  Y = V^2$. It is a radicial extension of degree 4. The norm
of the generic element
$T_0+T_1U+T_2V+T_3UV$ is \;$(T_0^2+T_1^2X+T_2^2Y+T_3^2XY)^2$. It
is a reducible polynomial.

\vspace{1cm}
\begin{center}
\textbf{3. Reading some pages of Hilbert.}
\end{center}
 \bigskip

The today reader of the beginning of the \emph{Zahlbericht} of Hilbert
(\cite{H}) has to face at least two difficulties. The first one is clearly
pointed out in the introduction
of the English edition; it comes from the deliberate avoiding by Hilbert of
any "abstract
algebra" concept, even the more useful ones among those he  already had at
hand (e.g the
notion of quotient group). The reader is thus driven through heavy
periphrases. In other
places, the constraint of rigour is less clear, for example when Hilbert
defines
something as the product of the (more or less mysterious) "conjugates" of some
expression, where today one understands simply a norm, etc.\par
The second difficulty is much more interesting because it is connected to
the debate,
at that time,  on the best way to "save" the factoriality (or, to say it
better, to get around the
non-factoriality) of the rings of integers (see the book by H. Weyl
\cite{W}, or the
Historical Note at the end of
\cite{AC}). By oversimplifying this debate we may personify the two positions
by Dedekind and Kronecker. The idea of Kronecker was to associate polynomials
to the
objects under consideration, in order to work inside polynomial rings over
$\mathbb{Z}$, or
$\mathbb{F}_p$, which are indeed factorial. The idea of Dedekind was to
create the
notion of ideal. It prevails.\bigskip

In what follows, I should like to show how the language and the results of
the preceding
paragraphs enlighten some of the statements of the
\emph{Zahlbericht}, mainly those which involve the Kronecker construction.
I think this
construction deserves to be better known even if it became useless in the
algebraic
theory of numbers.\par
The statements in question are in the \S\S 10 and 11; here the base ring is
$A = \mathbb{Z}$ and the algebra $B$ is the ring of integers of a number
field $K$. The generic element $\xi$ is called by Hilbert the \emph{fundamental
form}, and  the
\gcp\; is denoted by $F$ (amazingly enough, it is named in the memoir as
"the left hand
side of the fundamental equation"). As shown before,
$B$ is locally simple over
$A$, hence the Kronecker morphism is injective, and it remains injective
modulo any prime
$p$; this is the content of theorem 34 of the
\emph{Zahlbericht} which says:\\

{\it The congruence of degree} $n$,\; $F(X) \equiv 0$ mod.$p$ {\it is the
congruence
of lowest degree which is satisfied modulo $p$ by the fundamental form
$\xi$ (i.e
by the generic element).}\bigskip

The first part of theorem 33 of Hilbert's memoir gives the
correspondence between
the factorization of \ $F(X)$\ mod. $p$\, (that is in $\mathbb{F}_p[T_1,
\dots,
T_n, X]$) and the (now) usual factorization of the ideal $pB$ as a product
of prime ideals in $B$. Namely:\\

{\it If $p$ factorizes in $B$ as $pB = \mathfrak{p}^e{\mathfrak{p'}}^{e'}
\cdots$ then $F$ decomposes modulo $p$ in the form}
$$
F \equiv \Pi^e {\Pi'}^{e'} \cdots  \quad \mbox{mod.} p,
$$
{\it where $\Pi, {\Pi'},\dots$ represent distinct polynomials which are
irreducible modulo $p$.}\\

The proof in the memoir is preceded by three lemmas we may nowadays easily
circumvent by some functoriality considerations. In fact, the factorization of
$pB$ gives a decomposition of
$\mathbb{F}_p$-algebras:
$$
B/pB = B/\mathfrak{p}^e \times B/\mathfrak{p'}^{e'}\times \dots
$$
Using the remarks 2.3.3 and 2.3.4 (and without mentioning the base field
$\mathbb{F}_p$ any more), we may write
$$
F = v(G^{e})v'({G'}^{e'}) \cdots
$$
where $G$ is the \gcp \; of $B/\mathfrak{p}$, and where $v, v' \ldots$
denote the injective morphisms between rings of parameters
$$v :
\mbox{Sym}((B/\mathfrak{p})\,\check{}) \rightarrow
\mbox{Sym}((B/pB)\,\check{})$$
associated with the quotients $B/\mathfrak{p}, B/\mathfrak{p'} \ldots$. Since
these morphisms may be seen as just "adding the variables corresponding
to a basis of the
kernel", they preserve the irreducibility of polynomials. But the \gcp \;
$G$ is irreducible in
$\mbox{Sym}((B/\mathfrak{p}\,\check{})[X]$ ({\bf 2.4}). Therefore its image
$v(G)$ is still irreducible in
$\mbox{Sym}((B/pB)\,\check{})[X]$. This is the required polynomial $\Pi$.
\bigskip

Theorem 35 of the \emph{Zahlbericht} is also, as Hilbert
pointed out, a consequence of the injectivity of the Kronecker morphism
:\medskip

{\it The content of the discriminant of
$F(X)$ is equal to the discriminant of $B$ (or of $K$).}\\

The discriminant of $F(X)$ is an element of the ring containing the
coefficients
of $F$, i.e. here
$\mbox{Sym}(B\,\check{})$. This ring is isomorphic to $\mathbb{Z}[T_1, \dots,
T_n]$, therefore that makes sense to look at the gcd of the coefficients of the
discriminant, i.e at its \emph{content} (Hilbert writes: {\it the greatest
numerical factor}). Let us also recall what the discriminant of an algebra
is. Let $S\rightarrow E$ be a finite morphism, locally free
of rank $n$. The discriminant of $E/S$ is the ideal of $S$ image of the
$S$-linear map
$$
d_{E/S} : (\wedge^nE)^{\otimes 2} \longrightarrow S,
$$
defined as the extension to the $n$-th exterior power of the bilinear map
$$
E \times E \longrightarrow S, \qquad (x, y) \mapsto \mbox{Tr}_{E/S}(xy).
$$
If $F(X) \in S[X]$ is a monic polynomial, the discriminant of the $S$-algebra
$S[X]/(F)$ is the ideal generated by the discriminant of the polynomial
$F$. In the
situation under consideration, it can be checked that the Kronecker morphism
$$
u : E := S[X]/(F) \longrightarrow S\otimes_AB
$$
is compatible with the traces (see e.g \cite{F}, lemma 4.3.1 ), namely :
$$
\mbox{Tr}_{E/S} = \mbox{Tr}_{S\otimes_AB/S} \circ u .
$$
Since  $\mbox{Tr}_{S\otimes_AB/S} = \mbox{Tr}_{B/A}\otimes \mbox{id}_S$, we
get
$$
d_{E/S} = (d_{B/A}\otimes \mbox{id}_S)\circ (\wedge^n u)^{\otimes 2}.
$$
The Kronecker morphism $u$ is universally injective ({\bf 2.4}). Therefore
$\wedge^n u$ is injective, and remains
injective modulo any prime
$p$ (\cite{A} III 8.2 Prop.3). Hence, the content of $\wedge^n u$ is 1, and
the assertion follows.

\medskip
\textsc{IRMAR, Universit\'e de Rennes 1, Campus de Beaulieu, F-35040 Rennes
Cedex}\\
{\it E-mail address}\texttt{:\, daniel.ferrand@univ-rennes1.fr}
\end{document}